\documentclass[11pt]{article}

\usepackage{graphicx}
\usepackage{amsmath}
\usepackage{amsfonts}
\usepackage{fancyhdr} 
\usepackage{mathrsfs}

\usepackage{wrapfig}

\usepackage{bbold}

\usepackage{enumitem}

\newtheorem{lem}{Lemma}
\newtheorem{thm}{Theorem}

\newtheorem{prop}[lem]{Proposition}

\newtheorem{defn}[lem]{Definition}

\numberwithin{equation}{section}
\numberwithin{lem}{section}
\numberwithin{alg}{section}

\DeclareMathOperator{\sinc}{sinc}

\newcommand{\real}{\mathbb{R}}
\newcommand{\integer}{\mathbb{Z}}

\newcommand{\ddd}{\mathbb{D}}

\newcommand{\complex}{\mathbb{C}}

\newenvironment{pf}{\noindent {\em Proof}.\ \ }{\hspace*{\fill}\rule{.5ex}{1.4ex}\,}

\newcommand{\lll}{\mathcal L}
\newcommand{\one}{\mathbb{1}}

\newcommand{\ccc}{\mathcal{C}}

\newcommand{\ttt}{\mathbb{T}}

\newcommand{\dist}{\mathcal{D}^\prime}

\newcommand{\ml}{\widetilde{\Delta}}

\DeclareMathOperator{\spn}{span}

\title{\vspace*{-1in}{F}ourier expansion of disk automorphisms via scattering in layered media}
\author{Peter C.~Gibson\footnote{Dept.~of Mathematics \& Statistics, York University, 4700 Keele St., Toronto, Ontario, Canada, M3J~1P3, $\mathtt{pcgibson@yorku.ca}$} }
\date{March 7, 2016}
\makeatletter
\let\newtitle\@title
\let\newauthor\@author
\let\newdate\@date
\makeatother

\begin{document}
\maketitle
\begin{abstract}
A family of orthogonal polynomials on the disk (which we call scattering polynomials) serves to formulate a remarkable Fourier expansion of the composition of a sequence of Poincar\'e disk automorphisms.  Scattering polynomials are tied to an exotic riemannian structure on the disk that is hybrid between hyperbolic and euclidean geometries, and the expansion therefore links this exotic structure to the usual hyperbolic one.  The resulting identity is intimately connected with the scattering of plane waves in piecewise constant layered media.  Indeed, a recently established combinatorial analysis of scattering sequences provides a key ingredient of the proof.  At the same time, the polynomial obtained by truncation of the Fourier expansion elegantly encodes the structure of the nonlinear measurement operator associated with the finite time duration scattering experiment.

\end{abstract}


\section{Introduction}

Using recent results from scattering theory we derive the Fourier series expansion for an $n$-fold composition of disk automorphisms, and show that the resulting identity connects diverse mathematical objects not previously known to be related.  In particular, scattering polynomials (an exotic family of orthogonal polynomials on the disk), a system of PDE with smooth coefficients, two different riemannian structures on the disk, and the scattering of acoustic waves in layered media are all shown to be connected through this identity.  

The main result is stated below in \S\ref{sec-main-result}.  In \S\ref{sec-layered} we briefly outline scattering in layered media and describe needed ingredients such as the reflection experiment Green's function.  We discuss related literature and prior art in \S\ref{sec-existing}.  The main result is proved in \S\ref{sec-proof}, some longer technical arguments being relegated to the appendices.  We explain in \S\ref{sec-structure} how the principal identity encodes the time limited scattering experiment and thereby provides fundamental insights into its geometric and computational aspects.  The paper concludes with a brief discussion in \S\ref{sec-conclusion}. 

\subsection{The main result\label{sec-main-result}}
Let $\integer_+$ denote the nonnegative integers, including 0, let $\ddd\subset\complex$ denote the unit disk in the complex plane, $\overline{\ddd}$ its closure, and $\ttt$ its boundary, the unit circle.  Given a pair of $n$-tuples $(w,z)\in\overline{\ddd}^n\times\ttt^n$, let $\Psi^{w_j}_{z_j}:\overline{\ddd}\rightarrow\overline{\ddd}$ denote the linear fractional transformations
\begin{equation}\label{disk-automorphism}
\Psi^{w_j}_{z_j}(v)=z_j\frac{v+w_j}{1+\overline{w}_jv}\quad\quad(1\leq j\leq n)
\end{equation}
which if $w\in\ddd^n$ are automorphisms of the Poincar\'e disk.  (Note that if $w_j\in\ttt$ is on the disk boundary then $\Psi^{w_j}_{z_j}(v)=z_jw_j$ collapses to a constant function.)  We fix notation for the polar form of the entries of $w$ and $z$ as
\begin{equation}\label{polar-form}
w_j=r_je^{i\theta_j}\quad\mbox{ and }\quad z_j=e^{i\xi_j}\quad\quad(1\leq j\leq n).
\end{equation}
Our central object of interest is the value at 0 of an $n$-fold composition of the transformations (\ref{disk-automorphism}),
\begin{equation}\label{n-fold}
\Psi(w,z)=\Psi^{w_1}_{z_1}\circ\cdots\circ\Psi^{w_n}_{z_n}(0),
\end{equation}
viewed as a function of the parameters $(w,z)$.  For each fixed $w$, (\ref{n-fold}) is a function on the $n$-torus $z\in\ttt^n$ and has a corresponding Fourier series with coefficients depending on $w$.  It turns out that these coefficients have a beautiful structure, expressible very simply in terms of a new family of orthogonal polynomials which we call scattering polynomials. 
\pagestyle{fancyplain}
\begin{defn}[scattering polynomials]\label{defn-varphi}
For each $(p,q)\in\integer^2$ write $\varphi^{(p,q)}:\complex\rightarrow\complex$ to denote the polynomial defined as follows. If $\min\{p,q\}\geq 1$ set
\begin{equation}\label{varphi}
\varphi^{(p,q)}(\zeta)=\textstyle\frac{(-1)^p}{q(p+q-1)!}\,\displaystyle(1-\zeta\bar{\zeta})\frac{\partial^{\,p+q}}{\partial\bar{\zeta}^p\partial\zeta^q}(1-\zeta\bar{\zeta})^{p+q-1}.
\end{equation}
If $\min\{p,q\}<0$ or $p=0<q$ set $\varphi^{(p,q)}=0$; and if $p\geq 0$ set $\varphi^{(p,0)}(\zeta)=\zeta^p$.  
\end{defn}
That each $\varphi^{(p,q)}$ is actually a polynomial follows directly from (\ref{varphi}), as does the fact that $pq\varphi^{(p,q)}(\zeta)=0$ if $\zeta\in\ttt$ (i.e. on $\ttt$ $\varphi^{(p,q)}(\zeta)\neq0$ only if $pq=0$).  Less obvious is that the scattering polynomials $\varphi^{(p,q)}$ comprise a complete set of eigenfunctions of the weighted laplacian $\ml$ defined as
\begin{equation}\label{operator}
\ml=\frac{1-x^2-y^2}{4}\Delta=(1-\zeta\bar{\zeta})\frac{\partial^2}{\partial\bar{\zeta}\partial\zeta}
\end{equation}
where $\zeta=x+iy$ and $\Delta=\partial_{xx}+\partial_{yy}$.  (This and other properties of scattering polynomials and the associated operator $\ml$ are detailed in \S\ref{sec-complexification}.)   Note that $\ml$ may be viewed as hybrid between the euclidean laplacian $\Delta$ and the hyperbolic laplacian $(1-x^2-y^2)^2\Delta/4$, the Laplace-Beltrami operator for the Poincar\'e disk.  

Definition~\ref{defn-varphi} is all we need in order to describe the Fourier series for (\ref{n-fold}).  With respect to the various $n$-tuples that occur throughout the present paper, we adopt the convention that indices beyond $n$ should be interpreted as zero, so that for example, if $k\in\integer^n$, then $k_{n+1}=0$.   
\begin{thm}\label{thm-main} The Fourier series of the $n$-fold composition of disk automorphisms (\ref{n-fold}) is
\begin{equation}\label{main}
\Psi(w,z)=\sum_{k\in\{1\}\times\integer_+^n}\biggl(\,\prod_{j=1}^n\varphi^{(k_j,k_{j+1})}(w_j)\biggr)z^k.
\end{equation}
\end{thm}
(Note that invoking the notation (\ref{polar-form}), where $\xi=(\xi_1,\ldots,\xi_n)\in\real^n$, recovers the conventional form for Fourier characters $e^{i\langle\xi,k\rangle}=z^k$.) 

 Although it comes from acoustic imaging, specifically the scattering of waves in layered media, Theorem~\ref{thm-main} is of interest in its own right.  For instance, the right-hand side of (\ref{main}) is composed of eigenfunctions of the Laplace-Beltrami operator (i.e., $\ml$) for the riemannian metric $g_{ij}=4\delta_{ij}/(1-x^2-y^2)$, while the left-hand side of (\ref{main}) is comprised of automorphisms of the hyperbolic metric $g_{ij}=4\delta_{ij}/(1-x^2-y^2)^2$.  Thus the identity (\ref{main}) has a surprising geometric aspect---relating different riemannian structures to one another.  

The proof of Theorem~\ref{thm-main} relies on three ingredients: a system of linear PDE with smooth coefficients that has a unique solution, and which the left-hand side of (\ref{main}) can be shown to satisfy; and two results from scattering in layered media---one of which is a classic backwards recurrence formula, and the second of which is a combinatorial analysis of scattering sequences.  In the next section we describe more precisely the relevant layered media context.  At this point a word about notation is in order.  In the context of layered media we will follow common usage and use $z$ to denote a real scalar---the depth coordinate---and write $\zeta$ for acoustic impedance, which is a positive scalar-valued function of depth.  This is not to be confused with the very different usage above.

\subsection{Layered media\label{sec-layered}}

We consider layered media whose structure varies in the depth direction $z$ only; coordinates are oriented such that $z$ increases downward into the medium from a reference plane at depth 0. Given initial conditions that depend only on depth, the evolution of particle velocity $u(z,t)$ in a layered medium is governed by the one dimensional wave equation
\begin{equation}\label{wave}
u_{tt}-\frac{1}{\zeta}(\zeta u_x)_x=0
\end{equation}
where $\zeta(x)>0$ denotes acoustic impedance, and $x$ is travel time depth
\begin{equation}\label{x}
x(z)=\int_{0}^z\frac{1}{c(z^\prime)}\,dz^\prime.
\end{equation}
Suppose further that $\zeta(x)=C_0$ is constant for $x<0$, and that $\lim_{x\rightarrow\infty}\zeta(x)=1$.  A basic problem in acoustic imaging is to determine $\zeta(x)$ for $x>0$ from scattering data.  That is, suppose a source waveform $W\in C^\infty_c(\real)$ is supported on the negative half line, and initially travels downward, corresponding to initial conditions
\begin{equation}\label{initial}
u(x,0)=W(x)\quad u_t(x,0)=-W^\prime(x). 
\end{equation}
The resulting echoes $d(t)=u(0,t)$ recorded at $x=0$ during a certain time interval $0<t<T$ comprise scattering data that carries information about $\zeta(x)$ for $0<x<\frac{T}{2}$.  The (forward) scattering problem is to determine $d$ given $W$ and $\zeta$; the inverse scattering problem is to determine $\zeta$ given $d$, or given $d$ and $W$.  

We focus specifically on piecewise constant impedance of the form
\begin{equation}\label{impedance}
\zeta(x)=C_0+\sum_{j=1}^nC_jH(x-x_j).
\end{equation}
Here $H(x)$ denotes the Heaviside function, and jump points are indexed according to their natural order
\begin{equation}\label{x}
x_0=0<x_1<\cdots<x_n<x_{n+1}=\infty.
\end{equation}
The number of layers $n\geq 2$ is assumed to be fixed; but the particular value of $n$ can be arbitrarily large.  

Choosing $x_j^\ast\in(x_{j-1},x_j)$ arbitrarily for $1\leq j\leq n+1$, set 
\begin{equation}\label{Rtau}
r_j=\frac{\zeta(x_j^\ast)-\zeta(x_{j+1}^\ast)}{\zeta(x_j^\ast)+\zeta(x_{j+1}^\ast)}\quad\mbox{ and }\quad \tau_j=2(x_j-x_{j-1})\quad\mbox{ for }\quad 1\leq j\leq n.
\end{equation}
Writing
\begin{equation}\label{parameters}
(\tau,r)=\bigl((\tau_1,\ldots,\tau_n),(r_1,\ldots,r_n)\bigr),
\end{equation}
one can easily reconstruct $\zeta$ in (\ref{impedance}) from the pair $(\tau,r)$, giving a correspondence 
\begin{equation}\label{correspondence}
\zeta\leftrightarrow(\tau,r).
\end{equation}
One may thus regard the reflectivity vector $r$ and the travel time per layer vector $\tau$ as physical parameters characterizing a given piecewise constant layered medium.  It turns out that the measured data $d$ for the scattering problem can be expressed very simply in terms of $(\tau,r)$; for convenience we therefore 
work with $(\tau,r)$ instead of $\zeta$.  

Standard analysis of reflection and transmission at layer interfaces shows that measured scattering data has the general form 
\begin{equation}\label{measured-data}
d=\chi_{(0,T)}W\ast G^{(\tau,r)},
\end{equation}
where $G^{(\tau,r)}$ denotes the impulse response, or Green's function, at the boundary.  The latter is a delta train, expressible in the form
\begin{equation}\label{Green-form}
G^{(\tau,r)}(t)=\sum_{i=1}^\infty a_i\delta(t-t_i),
\end{equation}
where, for each $i$, $0<t_i<t_{i+1}$ and $a_i\neq 0$.  We refer to the points $t_i$ as arrival times (always indexed according to their natural order), and the corresponding $a_i$ as amplitudes.  

For present purposes we need only consider the impulsive limit $W\rightarrow\delta$, in which the measured data 
\begin{equation}\label{impulsive-limit-data}
d(t)=\chi_{(0,T)}G^{(\tau,r)}(t)=\sum_{i=1}^{i_{\max}}a_i\delta(t-t_i)\quad\quad i_{\max}=\max\{i\,|\,0<t_i<T\}
\end{equation}
is a time limited truncation of the Green's function $G^{(\tau,r)}$ itself.  

The objective of acoustic imaging is to solve the inverse problem, that is, to infer the physical parameters $(\tau,r)$ given the data $d$ in (\ref{impulsive-limit-data}).  This requires knowing how the amplitudes and arrival times depend on the physical parameters.  While the structure of the arrival times is reasonably straightforward (each $t_i$ has the form $\langle \tau,k\rangle$ for some $k\in\integer_+^n$), the amplitudes depend on the cumulative weight of successive reflections and transmissions from all simultaneously arriving scattering sequences.  The associated combinatorial problem was solved in \cite{Gi:SIAP2014}, relatively recently considering the history of the subject. (See \cite{Gi:SIAP2014} for further details on scattering sequences.)  The resulting formulas for the $a_i$ serve as a first step toward proving Theorem~\ref{thm-main}.  

\subsection{Connections to existing literature\label{sec-existing}}

Scattering polynomials $\varphi^{(p,q)}$ for which $\min\{p,q\}\geq 1$ are orthogonal with respect to the area measure $d\mu=4\,dxdy/(1-x^2-y^2)$ on the disk, but they do not belong to any of the classical families of orthogonal polynomials.  This is because monomials $x^my^n$ are not in general integrable with respect to $d\mu$ (even 1 is not integrable!), and so scattering polynomials cannot be generated according to the classical prescription of applying Gram-Schmidt orthogonalization to monomials.  Nevertheless, recent developments in approximation theory (independent from the work leading to the present paper) have produced general classes of orthogonal polynomials of which scattering polynomials can be see to be a special case; see \cite[Proposition~2.1 and Lemma~3.2]{LiXu:2014} and \cite[\S3]{CaKo:2015}.  The latter work does not contain the Rodrigues-type formula in Definition~\ref{defn-varphi} or the weighted laplacian $\ml$, nor does it mention the corresponding riemannian structure---so there is limited overlap.  The connection of scattering polynomials to scattering theory and to disk automorphisms presented here is completely new.  

There is an extensive literature on the one dimensional scattering problem associated with layered media; see \cite{BuBu:1983} for numerous early references. There are two main classes of media that have been considered, depending on the structure of the impedance function $\zeta$ occurring in (\ref{wave}).  The first class corresponds to  square integrable $\zeta^\prime/\zeta\in L^2(\real)$. This case is well understood, as exemplified by the paper \cite{SyWi:1999}.  Another class of media disjoint from the first corresponds to piecewise constant $\zeta$---the class we consider here.  Superficially, piecewise constant $\zeta$ seems like a much simpler structure than arbitrary $\zeta^\prime/\zeta\in L^2(\real)$, but in fact from the deterministic point of view piecewise constant media are less well understood.  Discontinuous coefficients in (\ref{wave}) have proven to be a formidable technical difficulty.  The main approach to piecewise constant media has been probabilistic, under the rubric of random media.  This has a large literature of its own, much of which is synthesized in the treatise \cite{FoGaPaSo:2007}.  The most recent development in deterministic analysis of piecewise constant layered media is \cite{Gi:SIAP2014}, which serves as a direct antecedent to the present paper.

\section{Proof of Theorem~\ref{thm-main}\label{sec-proof}}

The proof of Theorem~\ref{thm-main} is organized as follows.
We link together two known results from the theory of scattering in layered media to obtain a special case of the theorem.  This is carried out in \S\ref{sec-preliminaries}, Theorem~\ref{thm-special-case} being a version of the main result valid only for real parameters $w\in[-1,1]^n$.   The central technical difficulty is the extension of Theorem~\ref{thm-special-case} to complex $w\in\overline{\ddd}^n$.  

Next in \S\ref{sec-complexification} we establish several facts about the scattering polynomials, connecting them to both their real counterparts and to the modified laplacian $\ml$. 

In \S\ref{sec-system} we introduce the technical key, which is a system of PDE with smooth coefficients, and we complete the proof.  We show that $\Psi(w,z)$ satisfies the PDE (Lemma~\ref{lem-left}) and that the PDE has a unique solution in the form of the required Fourier expansion, up to a sequence of undetermined scalars (Lemma~\ref{lem-unique}).  Theorem~\ref{thm-main} follows once the scalars are determined, which is accomplished using Theorem~\ref{thm-special-case}.  The proofs of the two main lemmas are deferred to the appendices. 

\subsection{Preliminaries\label{sec-preliminaries}}

The purpose of this section is two recall two known results which combine to yield a special case of Theorem~\ref{thm-main}.   

 Here we use the (physicists') version of the Fourier transform consistent with the formula
\begin{equation}\label{Fourier}
\hat{f}(\sigma)=\int_{\real}f(t)e^{i\sigma t}\,dt.  
\end{equation}
The following well-know result is sometimes referred to as a backward recurrence formula.
\begin{thm}[See {\cite[\S3.5.2]{FoGaPaSo:2007}}]\label{thm-backward}
The Fourier transform with respect to time of the Green's function (\ref{Green-form}) is 
\[
\widehat{G^{(\tau,r)}}(\sigma)=\Psi_{e^{i\tau_1\sigma}}^{r_1}\circ\Psi_{e^{i\tau_2\sigma}}^{r_2}\circ\cdots\circ\Psi_{e^{i\tau_n\sigma}}^{r_n}(0).
\]
\end{thm}

The next result provides a key step leading to Theorem~\ref{thm-main}; it involves a family of polynomials of a real variable which were determined by combinatorial analysis and which led to the discovery of scattering polynomials. 
Given $(p,q)\in\integer^2$, set $m=|p-q|$ and $\nu=\min\{p,q\}-1$.  If $\min\{p,q\}\geq 1$, set 
\begin{equation}\label{fpq}
f^{(p,q)}(x)=\frac{(-1)^{q+\nu+1}}{q}(1-x^2)x^m\sum_{j=0}^{\nu}(-1)^j\frac{(j+\nu+m+1)!}{j!(j+m)!(\nu-j)!}x^{2j}
\end{equation}
and if $p\geq 0$ set $f^{(p,0)}(x)=x^p$; otherwise if $\min\{p,q\}<0$ or $p=0<q$, set $f^{(p,q)}(x)=0$.  
\begin{thm}[From {\cite[Theorems~2.4, 4.3]{Gi:SIAP2014}}]\label{thm-SIAP} The Green's function (\ref{Green-form}) has the form
\begin{equation}\label{Greens-formula}
G^{(\tau,r)}(t)=\sum_{k\in\{1\}\times\integer_+^{n-1}}\left(\prod_{j=1}^nf^{(k_j,k_{j+1})}(r_j)\right)\delta(t-\langle\tau,k \rangle).
\end{equation}
\end{thm}
Taking the Fourier transform of (\ref{Greens-formula}) with respect to $t$ yields the following special case of Theorem~\ref{thm-main}. 
\begin{thm}\label{thm-special-case}
For $(r,z)\in[-1,1]^n\times\ttt^n$ the composition (\ref{n-fold}) may be expanded as\begin{equation}\label{special-case}
\Psi(r,z)=\sum_{k\in\{1\}\times\integer_+^{n-1}}\left(\prod_{j=1}^nf^{(k_j,k_{j+1})}(r_j)\right)z^k.
\end{equation}
\end{thm}
\begin{pf}
Comparing Theorem~\ref{thm-backward} to the Fourier transform with respect to $t$ of Theorem~\ref{thm-SIAP} yields that for every $(\tau,r)\in\real_+^n\times[-1,1]^n$,
\begin{equation}\label{special-case-prelim}
\Psi_{e^{i\tau_1\sigma}}^{r_1}\circ\Psi_{e^{i\tau_2\sigma}}^{r_2}\circ\cdots\circ\Psi_{e^{i\tau_n\sigma}}^{r_n}(0)=\sum_{k\in\{1\}\times\integer_+^{n-1}}\left(\prod_{j=1}^nf^{(k_j,k_{j+1})}(r_j)\right)e^{i\langle \tau,k\rangle\sigma}.
\end{equation}
Since the mapping $(\tau,\sigma)\mapsto(e^{i\tau_1\sigma},\ldots,e^{i\tau_n\sigma})$ is a surjection from $\real_+^n\times\real$ onto $\ttt^n$, (\ref{special-case}) follows immediately from (\ref{special-case-prelim}).  
\end{pf}

There are considerable technical difficulties involved in progressing beyond Theorem~\ref{thm-special-case} to Theorem~\ref{thm-main}.  The first step is to recognize the correct complexification.  

\subsection{Complexification \label{sec-complexification}}

If we complexify the left-hand side of Theorem~\ref{thm-special-case} by extending the real parameters $r_j\in[-1,1]$ to complex parameters $w_j\in\overline{\ddd}$, then what is the corresponding modification of the right-hand side?  It clearly doesn't work to replace each $f^{(k_j,k_{j+1})}(r_j)$ with $f^{(k_j,k_{j+1})}(w_j)$, since the latter are holomorphic in $w_j$ whereas the disk automorphisms $\Psi^{w_j}_{z_j}$ depend on both $w_j$ and $\overline{w}_j$.    

Something more subtle is required.  What works is to replace each polynomial $f^{(k_j,k_{j+1})}(r_j)$ with $e^{i(k_j-k_{j+1})\theta_j}f^{(k_j,k_{j+1})}(r_j)$, which amounts to replacing each $f^{(k_j,k_{j+1})}(r_j)$ with the scattering polynomial $\varphi^{(k_j,k_{j+1})}(w_j)$, as we shall now show.  This is also a convenient juncture at which to verify the earlier claim that scattering polynomials are eigenvalues of the operator $\ml$ defined by the formulas (\ref{operator}).  
\begin{lem}\label{lem-complexification}
Referring to Definition~\ref{defn-varphi} and (\ref{fpq}), \begin{enumerate}
\item
for every $(p,q)\in\integer^2$ and $\rho e^{i\sigma}\in\complex$,
$\varphi^{(p,q)}(\rho e^{i\sigma})=e^{i(p-q)\sigma}f^{(p,q)}(\rho)$, 
and
\item 
for every $(p,q)\in\integer_+^2$, $-\ml\varphi^{(p,q)}=pq\,\varphi^{(p,q)}$.
\end{enumerate}
\end{lem}
\begin{pf}
Expanding the binomial $(1-\zeta\bar{\zeta})^{p+q-1}$ in the formula (\ref{varphi}), and then applying the derivative $\rule{0pt}{13pt}\partial^{p+q}/\partial\bar{\zeta}^p\partial \zeta^q$, yields
\begin{equation}\label{step-one}
\varphi^{(p,q)}(\zeta)=
\frac{(-1)^{q+\nu+1}}{q}(1-\zeta\bar{\zeta})\zeta^{m+\nu-q+1}\bar{\zeta}^{m+\nu-p+1}\sum_{j=0}^{\nu}(-1)^j\frac{(j+\nu+m+1)!}{j!(j+m)!(\nu-j)!}(\zeta\bar{\zeta})^j,
\end{equation}
where $m=|p-q|$ and $\nu=\min\{p,q\}-1$.  Switching to polar form $\zeta=\rho e^{i\sigma}$, it follows from (\ref{step-one}) that 
\begin{equation}\label{step-two}
\varphi^{(p,q)}\bigl(\rho e^{i\sigma}\bigr)=e^{i(p-q)\sigma}\frac{(-1)^{q+\nu+1}}{q}(1-\rho ^2)\rho ^m\sum_{j=0}^{\nu}(-1)^j\frac{(j+\nu+m+1)!}{j!(j+m)!(\nu-j)!}\rho ^{2j}.
\end{equation}
Comparing this with (\ref{fpq}) proves part 1. 

One can verify directly that $-\ml\varphi^{(p,q)}=pq\,\varphi^{(p,q)}$ by applying the operator 
\[
-\ml=-(1-\zeta\bar{\zeta})\frac{\partial^2}{\partial\zeta\partial\bar{\zeta}}
\]
to the right-hand side of (\ref{step-one}) to prove part 2.
\end{pf}

Much more is true.  The unit disk endowed with the riemannian metric $g_{ij}=4\delta_{ij}/(1-x^2-y^2)$, with respect to which $\ml$ is the Laplace-Beltrami operator, has a very rich structure.  But this is a topic for a separate paper.  For present purposes we confine ourselves to the following fact, whose proof is detailed in appendix~\ref{sec-app-hybrid}.  
\begin{thm}\label{thm-hybrid-laplacian}
Given Dirichlet boundary conditions on $\overline{\ddd}$, the non-zero eigenvalues of $-\ml$ are positive integers.  For each positive integer $\lambda$,
\[
\ker(\ml+\lambda)=\spn\bigl\{\varphi^{(p,q)}\,|\;pq=\lambda\mbox{ and }p,q\in\integer_+\bigr\}.
\]
In particular, the dimension of $\ker(\ml+\lambda)$ is the number of divisors of $\lambda$.  
\end{thm}
In the next section we need to consider the weighted laplacian with respect to each of the complex variables $w_j$, which we denote
\begin{equation}\label{mlj}
\ml_j=(1-w_j\overline{w}_j)\frac{\partial^2}{\partial\overline{w}_j\partial w_j}.
\end{equation}

\subsection{A smooth PDE; conclusion of the proof\label{sec-system}}

At this point it is convenient to streamline notation for products of scattering polynomials as follows.  Set
\begin{equation}\label{notation}
\Phi(w,k)=\prod_{j=1}^n\varphi^{(k_j,k_{j+1})}(w_j)\quad\quad\;\;(w,k)\in\overline{\ddd}^n\times\integer^n.
\end{equation}
Recall the notation (\ref{polar-form}) for the polar forms of $(w,z)\in\overline{\ddd}^n\times\ttt^n$.  Let $u:\overline{\ddd}^n\times\ttt^n\rightarrow\complex$ denote an unknown function.  The following system of equations plays a key technical role in proving Theorem~\ref{thm-main}.  Here $1\leq j\leq n$,  $\one=(1,1,\ldots,1)$ and $\partial/\partial\xi_{n+1}=0$.
\begin{subequations}\label{sys}
\begin{align}
-\widetilde{\Delta}_ju+\frac{\partial^2u}{\partial\xi_{j+1}\partial\xi_j}&=0\label{proof-wave}\\
\frac{\partial u}{\partial\theta_j}-\frac{\partial u}{\partial\xi_j}+\frac{\partial u}{\partial\xi_{j+1}}&=0\label{proof-coupling}\rule{0pt}{22pt}\\
\mbox{ if }w_n\in\ttt \mbox{ then }\; u(w,\one)&=\Psi(w,\one)\label{proof-boundary}\rule{0pt}{22pt}
\end{align}
\end{subequations}

\begin{lem}\label{lem-left}
The $n$-fold composition of disk automorphisms $\Psi(w,z)$ satisfies (\ref{sys}).  
\end{lem}
Verification of Lemma~\ref{lem-left} is somewhat technical and so is consigned to appendix~\ref{sec-app-left}.   The next result is that the system (\ref{sys}) has a unique solution that can be computed in series form directly by separation of variables.  This is possible precisely because the equations (\ref{sys}) have real analytic coefficients---in marked contrast to the discontinuous coefficient that appears in the wave equation (\ref{wave}).  Separation of variables leads to an expansion of the solution to (\ref{sys}) in terms of eigenfunctions of the operators $\ml_j$, that is, in terms of scattering polynomials up to a sequence of scalars.  The scalars can then be determined using the embryonic form, Theorem~\ref{thm-special-case}, of the main result.  
\begin{lem}\label{lem-unique}
The system (\ref{sys}) has a unique distributional solution, necessarily of the form
\begin{equation}\label{form}
u(w,z)=\sum_{k\in\{1\}\times\integer_+^n}\beta_k\Phi(w,k)\,z^k
\end{equation}
where each $\beta_k$ is a scalar.  
\end{lem}
The proof of Lemma~\ref{lem-unique} is also rather technical and is deferred to appendix~\ref{sec-app-unique}.  

Lemmas~\ref{lem-left} and \ref{lem-unique} combine to show that 
\begin{equation}\label{nearly}
\Psi(w,z)=\sum_{k\in\{1\}\times\integer_+^n}\beta_k\Phi(w,k)z^k
\end{equation}
for some collection of scalars $\beta_k$.  It follows from Theorem~\ref{thm-special-case} and the first part of Lemma~\ref{lem-complexification} that for every $r\in[-1,1]^n$, 
\begin{equation}\label{reduces}
\sum_{k\in\{1\}\times\integer_+^n}\Phi(r,k)z^k=\sum_{k\in\{1\}\times\integer_+^n}\beta_k\Phi(r,k)z^k.
\end{equation}
Since for each $k\in\integer^n$, $\Phi(r,k)\not\equiv0$ if $\Phi(w,k)\not\equiv0$, equation (\ref{reduces}) implies by uniqueness of Fourier coefficients that $\beta_k=1$ for each $k$ such that $\Phi(w,k)\not\equiv0$, thereby reducing (\ref{nearly}) to Theorem~\ref{thm-main}.

\section{Structure of the scattering data\label{sec-structure}}

In the present section we return to scattering in layered media and discuss the connection of the time-limited truncation of the Green's function, $\chi_{(0,T)}G^{(\tau,r)}$, to the Fourier expansion of $\Psi(w,z)$ given in Theorem~\ref{thm-main}.  The emphasis is primarily on qualitative insights. 

To begin, we compare previously known results, in the guise of Theorem~\ref{thm-backward}, to Theorem~\ref{thm-main} and its precursor Theorem~\ref{thm-special-case}.  
Translating Theorem~\ref{thm-backward} into the current the current notation (\ref{notation}) yields the formula
\begin{equation}\label{translation-to-current}
\widehat{G^{(\tau,r)}}(\sigma)=\Psi(r,z(\sigma))\quad\mbox{ where }\quad z(\sigma)=\bigl(e^{i\tau_1\sigma},\ldots,e^{i\tau_n\sigma}\bigr). 
\end{equation}
In order to use this formula to compute the Fourier transform of data measured over a finite time interval $0<t<T$ one has to contend with the convolution of the backward recurrence formula with a modulated sinc function, 
\begin{equation}\label{convolution}
\widehat{\chi_{(0,T)}G^{(\tau,r)}}(\sigma)=\frac{T}{2\pi}\int_{-\infty}^\infty \Psi(r,z(s))e^{i(\sigma-s)T/2}\sinc((\sigma-s)T/2)\,ds.
\end{equation}
Evaluation of this integral for a single $\sigma$ involves the values of $\Psi(r,z(s))$ for all $-\infty<s<\infty$, which hugely inflates the computational expense. 

By contrast, Theorem~\ref{thm-SIAP} in conjunction with Lemma~\ref{lem-complexification} and the notation (\ref{notation}) implies that time limited data may be expressed as 
\begin{equation}\label{time-limited-expression}
\chi_{(0,T)}G^{(\tau,r)}(t)=\sum_{\substack{k\in\{1\}\times\integer_+^{n-1}\\
\langle \tau,k\rangle<T}}\Phi(r,k)\delta(t-\langle\tau,k \rangle),
\end{equation}
whose Fourier transform is 
\begin{equation}\label{time-limited-fourier}
\widehat{\chi_{(0,T)}G^{(\tau,r)}}(\sigma)=\sum_{\substack{k\in\{1\}\times\integer_+^{n-1}\\
\langle \tau,k\rangle<T}}\Phi(r,k)z(\sigma)^k.  
\end{equation}
There is a finite set of lattice points $k\in\integer_+^n$ that meet the constraint $\langle \tau,k\rangle<T$, so the right-hand side of (\ref{time-limited-fourier}) is a polynomial.  This shows that (\ref{convolution})
may be computed \emph{exactly} by polynomial evaluation, avoiding altogether the necessity for approximation imposed by the convolution (\ref{convolution}). (In addition to being exact, polynomial evaluation is inherently fast for low dimension $n$ and sufficiently short times $T$.)  

Equation (\ref{time-limited-fourier}) shows how Theorem~\ref{thm-main} encodes the time limited scattering experiment: truncation of the Fourier series for $\Psi(w,z)$, evaluated at $(r,z(\sigma))$, is the Fourier transform of the time limited scattering data.  Indeed this encoding illuminates the previously hidden geometric nature of the scattering process, as follows.  Note that the function 
\begin{equation}\label{line}
z(\sigma)=\bigl(e^{i\tau_1\sigma},\ldots,e^{i\tau_n\sigma}\bigr)\quad\quad-\infty<\sigma<\infty
\end{equation}
is a line on the torus with direction vector $\tau$.  Generically, in particular if $\tau_1,\ldots,\tau_n$ are linearly independent over $\integer$, the line (\ref{line}) is a leaf of the $n$-dimensional analogue of a Kronecker foliation, and therefore everywhere dense on $\ttt^n$.  The scattering data samples the polynomial 
\begin{equation}\label{polynomial}
P_T(r,z)=\sum_{\substack{k\in\{1\}\times\integer_+^{n-1}\\
\langle \tau,k\rangle<T}}\Phi(r,k)z^k
\end{equation}
on the leaf $z(\sigma)$, which winds around the torus haphazardly, according to the ergodic nature of the generalized Kronecker foliation, eventually sampling $P_T(r,z)$ near every $z\in\ttt^n$.  (See \cite[\S1.5 and \S4.2]{KaHa:1995}.)   This remarkable structure is completely hidden in the scattering data in its raw form as described by equation (\ref{impulsive-limit-data}), which just looks like noise.  

The scattering data fails to sample $P_T(r,z)$ densely on the torus in the non-generic case where the entries of $\tau$ are all integer multiples of a fixed quantity.  From the point of view of trying to recover $r$ from the scattering data this is the worst possible scenario---yet it corresponds to Gaupillaud media, perhaps the most extensively studied model for piecewise constant layered media (see \cite[\S3.5.4]{FoGaPaSo:2007} and \cite{BuBu:1983}).  Gaupillaud media have layers of equal thickness.  This seems like a natural simplification from the point of view of ultimate discretization (the rationale in \cite{BuBu:1983}) but the geometric discussion above shows that this simplification misses the generic behaviour of the underlying model.  

In summary, Theorem~\ref{thm-main} encodes the scattering of waves in layered media in a way that gives substantial insight into the finite time duration scattering experiment both from a computational and geometric perspective---and exposes a drawback of the common assumption of constant layer depth.     

\section{Conclusion\label{sec-conclusion}}

In the present section we discuss very briefly the contributions of the paper and suggest some directions for further investigation. 

A principal contribution of the present paper, beyond Theorem~\ref{thm-main} itself, is to show that the system of PDE (\ref{sys}) has a direct bearing on the propagation of waves in layered media. The weighted laplacian $\ml$ occurring in (\ref{sys}) and its eigenfunctions, the scattering polynomials, were not previously known to have any such role. In fact, as discussed in \S\ref{sec-existing}, scattering polynomials themselves are new, having only recently appeared in the approximation theory literature in the context of general families of orthogonal polynomials. And the riemannian structure on the disk $g_{ij}=4\delta_{ij}/(1-x^2-y^2)$, of which $\ml$ is the Laplace-Beltrami operator, seems to be entirely unknown.  

From the  perspective of PDE generally, there is relatively little machinery available for handling discontinuous coefficients.  We have demonstrated that the discontinuous coefficient in (\ref{wave}) can be handled by means of the smooth higher dimensional system (\ref{sys}), in which there is one dimension for each discontinuity in the original coefficient.  This is highly suggestive: can the idea of trading a PDE with discontinuous coefficients for a smooth higher dimensional system be implemented in other contexts? 

Theorem~\ref{thm-main} is also interesting quite apart from its relation to layered media.  From the geometric point of view, it is surprising that automorphisms of the Poincar\'e disk should have such a seemingly natural connection to orthogonal functions on the disk with metric tensor $g_{ij}=4\delta_{ij}/(1-x^2-y^2)$---in other words, scattering polynomials.  It would be interesting to know whether this is a special case of a more general phenomenon.

\appendix

\section{Some technical proofs\label{sec-technical}}

\subsection{Proof of Theorem~\ref{thm-hybrid-laplacian}\label{sec-app-hybrid}}

Note that the angular part of $\varphi^{(p,q)}(\rho e^{i\sigma})$, namely $e^{i(p-q)\sigma}$, is a pure frequency.  Therefore if $p-q\neq p^\prime-q^\prime$, then $\varphi^{(p,q)}$ and $\varphi^{(p^\prime,q^\prime)}$ are orthogonal, both in $L^2\bigl(\ddd,4dxdy/(1-x^2-y^2)\bigr)$ and $L^2\bigl(\ddd,dxdy\bigr)$.  In particular, if $pq=p^\prime q^\prime$ and $(p,q)\neq(p^\prime,q^\prime)$, then $\varphi^{(p,q)}$ and $\varphi^{(p^\prime,q^\prime)}$ are orthogonal, so the set of scattering polynomials corresponding to any fixed eigenvalue is linearly independent.  

In order to complete the proof of Theorem~\ref{thm-hybrid-laplacian}, it remains to show that: (1) only non-negative integers are eigenvalues of $-\ml$ with zero boundary values; and (2) for any two integers $n\neq0$ and $\lambda\geq 1$ there is at most one radial function $f(\rho )$ (up to scalar multiplication) such that $f(\rho )e^{in\sigma}$ is an eigenfunction of $-\ml$ with eigenvalue $\lambda$.  Both (1) and (2) will be seen to follow from separation of variables applied to the eigenvalue equation for $-\ml$, as follows.   

Since it is elliptic and has analytic coefficients, the operator $\ml+\lambda$ is analytic hypoelliptic for any constant $\lambda$. Therefore any distributional eigenfunction $\varphi$ of $-\ml$ is necessarily a real analytic function \cite[Thm.~10]{Ha:2006}.  Such a function is the uniform limit of its radial Fourier series
\[
\varphi(\rho e^{i\sigma})=\sum_{n\in\integer}a_n(\rho )e^{in\sigma},
\]
on which the operator $-\ml$ may be evaluated term by term.   The image of any particular term $a_n(\rho )e^{in\sigma}$ by $-\ml$ is a function $A_n(\rho )e^{in\sigma}$ having the same angular part.  Hence if $\varphi$ is an eigenfunction of $-\ml$ with eigenvalue $\lambda$, the same is true of each nonzero term $a_n(\rho )e^{in\sigma}$.  In other words, real analytic tensor products of the form $f(\rho )e^{in\sigma}$, where $n\in\integer$, span the eigenspaces of $-\ml$, and separation of variables is guaranteed not to miss any solutions to the eigenvalue problem.  

Suppose therefore that, for some $n\in\integer$, $f(\rho )e^{in\sigma}$ is a real analytic eigenfunction of $-\ml$ corresponding to non-zero eigenvalue $\lambda$ such that $f(1)=0$.  We shall verify that the eigenvalue equation determines $f$ up to a scalar multiple, and that $\lambda$ is necessarily a positive integer, as follows.  Expressing the hybrid laplacian in polar coordinates yields
\[
(\ml+\lambda)\bigl(f(\rho )e^{in\sigma}\bigr)=\frac{1-\rho^2}{4}\left(f^{\prime\prime}(\rho )+\frac{1}{\rho}f^\prime(\rho )+\left(\frac{4\lambda}{1-\rho^2}-\frac{n^2}{\rho^2}\right)f(\rho )\right)e^{in\sigma}.
\]
Therefore $(\ml+\lambda)\bigl(f(\rho )e^{in\sigma}\bigr)=0$ implies that $f(\rho )$ satisfies the equation
\begin{equation}\label{ode}
\rho^2(1-\rho^2)f^{\prime\prime}+r(1-\rho^2)f^\prime+\bigl(4\lambda \rho^2-n^2(1-\rho^2)\bigr)f=0.
\end{equation}
By real analyticity $f$ has a convergent Taylor expansion,
\begin{equation}\label{taylor-expansion}
f(\rho )=\sum_{j=0}^\infty b_j\rho^j
\end{equation}
upon which equation (\ref{ode}) induces the recurrence relation
\begin{align}
n^2b_0&=0\label{b0}\\
(1-n^2)b_1&=0\label{b1}\\
\forall j\geq0\quad\bigl((j+2)^2-n^2\bigr)b_{j+2}&=\bigl(j^2-4\lambda -n^2)b_j\label{bj+2}.
\end{align}
Assuming $f$ is not identically zero, one deduces directly from (\ref{b0}),(\ref{b1}),(\ref{bj+2}) that the least index $m$ for which $b_m\neq 0$ is $m=|n|$, with the particular choice of $b_{|n|}$ determining $f$ itself by recurrence.  This proves that there is at most one radial function $f(\rho )$ (up to choice of $b_{|n|}$) such that $f(\rho )e^{in\sigma}$ is an eigenfunction of $-\ml$ corresponding to a given non-zero eigenvalue $\lambda $.   

Write $m=|n|$. Solving the recurrence yields $b_{m+2j+1}=0$ $(j\geq 0)$, 
\begin{equation}\label{recurrence-solution}
b_{m+2}=\frac{-\lambda }{m+1}b_m\quad\mbox{ and }\quad 
b_{m+2j}=\frac{-\lambda }{j(m+j)}\prod_{\nu=1}^{j-1}\left(1-\frac{\lambda }{\nu(m+\nu)}\right)b_m\quad\quad\forall j\geq2.
\end{equation}
The boundary value $f(1)=0$ determines possible values of $\lambda $.  By (\ref{recurrence-solution}), $f(1)=b_m\xi(\lambda )$, where 
\begin{equation}\label{xi}
\begin{split}
\xi(\lambda )&=1-\lambda \left(\frac{1}{m+1}+\sum_{j=2}^\infty\frac{1}{j(m+j)}\prod_{\nu=1}^{j-1}\left(1-\frac{\lambda }{\nu(m+\nu)}\right)\right)\\
&=\prod_{\nu=1}^\infty\left(1-\frac{\lambda }{\nu(m+\nu)}\right).
\end{split}
\end{equation}
Convergence of the series $\sum1/(\nu(m+\nu))$ for all $\nu\geq1$ ensures that $f(1)=b_m\xi(\lambda )=0$ if and only if $\lambda =\nu(m+\nu)$ for some $\nu\geq1$.  This proves that $\lambda $ is a positive integer, completing the proof of Theorem~\ref{thm-hybrid-laplacian}.  

Slightly more work produces the scattering polynomials themselves, up to a scalar multiple.  It follows from (\ref{recurrence-solution}) that 
\[
b_{m+2j}=0\mbox{ if }j>\nu\quad\mbox{ and }\quad b_{m+2j}\neq0\mbox{ if }0\leq j\leq\nu. 
\]
Thus $f(\rho )$ is a polynomial of precise degree $m+2\nu$ that has a zero of order $m$ at $\rho=0$, and such a solution exists for every pair of integers $m\geq0$ and $\nu\geq 1$.  Set
\begin{equation}\label{beta}
\beta_0=1,\qquad\beta_j=1-\frac{\nu(m+\nu)}{j(m+j)}\quad(1\leq j\leq \nu).
\end{equation}
Note in particular that $\beta_\nu=0$.   Combining the formula (\ref{recurrence-solution}) with the expansion (\ref{taylor-expansion}) shows that 
\begin{equation}\label{explicit}
f(\rho )=b_m\rho^m(1-\rho^2)\sum_{j=0}^{\nu-1}\Bigl(\prod_{s=0}^j\beta_s\Bigr)\rho^{2j},
\end{equation}
which has two associated eigenvalues $\varphi(\rho e^{\pm im\sigma})=f(\rho )e^{\pm im\sigma}$ if $m\geq1$, and just one if $m=0$. With $m=|n|=|p-q|$ and $\nu=\min\{p,q\}$, the formula (\ref{explicit}) is proportional to $f^{(p,q)}(\rho )$ and $f^{(q,p)}(\rho )$ as defined in (\ref{fpq}).  Thus separation of variables yields the scattering polynomials (up to a scalar multiple) directly from the eigenvalue equation for $-\ml$.

\subsection{Proof of Lemma~\ref{lem-left}\label{sec-app-left}}

We shall verify that $\Psi(w,z) $ is a solution to (\ref{sys}), fixing notation for the relevant operators as
\[
\begin{split}
L_j&=-\widetilde{\Delta}_j+\frac{\partial^2}{\partial\xi_{j+1}\partial\xi_j}\\
C_j&=\frac{\partial }{\partial\theta_j}-\frac{\partial }{\partial\xi_j}+\frac{\partial }{\partial\xi_{j+1}}
\end{split}
\]
where $1\leq j\leq n$.

To begin we introduce a nonlinear first-order operator that will play an important auxiliary role in our analysis.  For each $1\leq j\leq n$, let $\nabla_j$ denote the gradient operator with respect to $w_j=x_j+iy_j\in\ddd$, so that 
\[
\nabla_jv=\left(\frac{\partial v}{\partial x_j},\frac{\partial v}{\partial y_j}\right)
\]
and define the differential operator $E_j$ on distributions on $\overline{\ddd}^n\times\ttt$ by the formula
\[
E_jv=-\frac{1-r_j^2}{4}\left(\nabla_jv\cdot\nabla_jv\right)+\frac{\partial v}{\partial\xi_{j+1}}\frac{\partial v}{\partial \xi_j}.
\]
(This is essentially the quadratic form associated with $L_j$, in the sense that $\langle L_jv,v\rangle=\int E_jv$.)   
\begin{prop}\label{prop-E_j}
$E_j\Psi(w,z) =0$ for every $1\leq j\leq n$.  
\end{prop}
\begin{pf}
We argue by induction on the sequence of functions
\[
v=\Psi^{w_j}_{z_j}\circ\cdots\circ\Psi^{w_n}_{z_n}(0)
\]
as $j$ decreases from $n$ to $1$.  The base case $v=\Psi^{w_n}_{z_n}(0)=z_nw_n$ trivially satisfies $E_jv=0$ if $j<n$.   Since $\frac{\partial v}{\partial\xi_{n+1}}=0$ (the variable $\xi_{n+1}$ does not occur), for $j=n$ it suffices to verify that 
\[
-\frac{1-r_n^2}{4}\left(\nabla_nv\cdot\nabla_nv\right)=0,
\]
which follows from the fact that $\nabla_nv=z_n(1,i)$.   

Next suppose that $1<j+1\leq n$ and that $E_sv^\prime=0$ for every $1\leq s\leq n$, where 
\[
v^\prime=\Psi^{w_{j+1}}_{z_{j+1}}\circ\cdots\circ\Psi^{w_n}_{z_n}(0).
\]
We will show that necessarily $E_sv=0$ where 
\[
v=\Psi^{w_j}_{z_j}\circ\cdots\circ\Psi^{w_n}_{z_n}(0).
\]
We may assume that $s\geq j$, since otherwise there is nothing to prove. 

 Suppose first that $s>j$.  By definition, 
\[
v=\Psi_{z_j}^{w_j}(v^\prime)=z_j\frac{w_j+v^\prime}{1+\overline{w}_jv^\prime}.
\]
It follows by direct computation that
\[
-\frac{1-r_s^2}{4}\left(\nabla_sv\cdot\nabla_sv\right)=
z_j^2\frac{(1-|w_j|^2)^2}{(1+\overline{w}_jv^\prime)^4}\frac{-(1-r_s^2)}{4}\left(\nabla_sv^\prime\cdot\nabla_sv^\prime\right),
\]
while
\[
\frac{\partial v}{\partial\xi_{s+1}}\frac{\partial v}{\partial\xi_{s}}=z_j^2\frac{(1-|w_j|^2)^2}{(1+\overline{w}_jv^\prime)^4}\frac{\partial v^\prime}{\partial\xi_{s+1}}\frac{\partial v^\prime}{\partial\xi_{s}}.
\]
The induction hypothesis therefore guarantees that $E_sv=0$.   

It remains to consider the case $s=j$.  Noting that $r_j=|w_j|$, straightforward computation yields that 
\[
-\frac{1-r_j^2}{4}\left(\nabla_jv\cdot\nabla_jv\right)=
\frac{(1-r_j^2)z_j^2(w_j+v^\prime)v^\prime}{(1+\overline{w}_jv^\prime)^3},
\]
while
\[
\frac{\partial v}{\partial\xi_{j+1}}\frac{\partial v}{\partial\xi_{j}}=\frac{-z_j^2(1-r_j^2)(w_j+v^\prime)v^\prime}{(1+\overline{w}_jjv^\prime)^3},
\]
from which follows the desired result that $E_jv=0$. 
\end{pf}

\begin{prop}\label{prop-L_j}
$L_j\Psi(w,z) =0$ for every $1\leq j\leq n$.
\end{prop}
\begin{pf}
This is proved along the same lines as Prop~\ref{prop-E_j}, by an induction on the sequence of functions
\[
v=\Psi^{w_j}_{z_j}\circ\cdots\circ\Psi^{w_n}_{z_n}(0)
\]
as $j$ decreases from $n$ to $1$.  The base case $v=\Psi^{w_n}_{z_n}(0)=z_nw_n$ trivially satisfies $L_jv=0$ if $j<n$.   Since $\frac{\partial v}{\partial\xi_{n+1}}=0$, for $j=n$ it suffices to verify that 
\[
-\widetilde{\Delta}_j=0,
\]
which follows from the fact that $z_nw_n$ is harmonic with respect to $w_n$.   

Next suppose that $1<j+1\leq n$ and that $L_sv^\prime=0$ for every $1\leq s\leq n$, where 
\[
v^\prime=\Psi^{w_{j+1}}_{z_{j+1}}\circ\cdots\circ\Psi^{w_n}_{z_n}(0).
\]
We will show that necessarily $L_sv=0$ where 
\[
v=\Psi^{w_j}_{z_j}\circ\cdots\circ\Psi^{w_n}_{z_n}(0).
\]
We may assume that $s\geq j$, since otherwise there is nothing to prove. 

 Suppose first that $s>j$.  Using that  
\[
v=\Psi_{z_j}^{w_j}(v^\prime)=z_j\frac{w_j+v^\prime}{1+\overline{w}_jv^\prime},
\]
direct computation shows
\[
-\widetilde{\Delta}_sv=z_j(1-|w_j|^2)\left(\frac{-\widetilde{\Delta}_sv^\prime}{(1+\overline{w}_jv^\prime)^2}-\frac{-2\overline{w}_j\frac{(1-r_s^2)}{4}\left(\nabla_sv^\prime\cdot\nabla_sv^\prime\right)}{(1+\overline{w}_jv^\prime)^3}\right),
\]
while
\[
\frac{\partial^2v}{\partial\xi_{s+1}\partial\xi_s}=z_j(1-|w_j|^2)\left(\frac{\frac{\partial^2v^\prime}{\partial\xi_{s+1}\partial\xi_s}}{(1+\overline{w}_jv^\prime)^2}-\frac{2\overline{w}_j\frac{\partial v^\prime}{\partial\xi_{s+1}}\frac{\partial v^\prime}{\partial\xi_{s}}}{(1+\overline{w}_jv^\prime)^3}\right).
\]
Adding the above two parts, the induction hypothesis then implies
\[
L_sv=\frac{-2z_j\overline{w}_j(1-|w_j|^2)}{(1+\overline{w}_jv^\prime)^3}E_sv^\prime,
\]
which is 0 by Prop~\ref{prop-E_j}.   

It remains to consider the case $s=j$.  Straightforward computation yields
\[
-\widetilde{\Delta}_jv=\frac{(1-|w_j|^2)z_jv^\prime}{(1+\overline{w}_jv^\prime)},
\]
while
\[
\frac{\partial^2v}{\partial\xi_{j+1}\partial\xi_j}=-\frac{(1-|w_j|^2)z_jv^\prime}{(1+\overline{w}_jv^\prime)^2},
\]
whence $L_jv=0$, completing the proof. \end{pf}

\begin{prop}\label{prop-C_j}
$\ccc_j\Psi(w,z) =0$ for every $1\leq j\leq n$.
\end{prop}
\begin{pf}
As in the foregoing results, this may be proved by downward induction on 
\[
v=\Psi^{w_j}_{z_j}\circ\cdots\circ\Psi^{w_n}_{z_n}(0).
\]
The appropriate induction hypothesis in this case is that $C_sv=0$ for every $s\geq j$.  It then follows by direct computation that $C_s\Psi_{z_{j-1}}^{w_{j-1}}(v)=0$ for every $s\geq j-1$, yielding the desired result.  (It is straightforward to check the only nontrivial case, $s=j-1$.)
\end{pf}

Lastly, $\Psi(w,z) $ satisfies the partial boundary condition (\ref{proof-boundary}) automatically and hence, by Propositions~\ref{prop-L_j} and \ref{prop-C_j}, the full system (\ref{sys}).

\subsection{Proof of Lemma~\ref{lem-unique}\label{sec-app-unique}}

The pairing of a distribution with a test function will be denoted by square brackets, whatever the domain; we take distributions to be continuous linear functionals (as opposed to conjugate linear).   Let $u\in\dist(\overline{\ddd}^n\times\ttt^n)$ be a distribution.    The Fourier coefficients of $u$ are distributions $c_k\in\dist(\overline{\ddd}^n)$ defined by the rule
\[
[c_k,\phi]=[u,\phi\otimes\mu_k]\qquad\bigl(k\in\integer^n,\;\phi\in C^\infty(\overline{\ddd}^n)\bigr),
\]
where $\mu_k(z)=\overline{z}^k/(2\pi)^n=e^{-i\langle k,\xi\rangle}/(2\pi)^n$.   
In order that $u$ be considered as a candidate solution to the system (\ref{sys}) it has to be possible to interpret the partial boundary condition (\ref{proof-boundary}). This requires in particular that $u$ be trace class in the sense that the series $\sum_{k\in\integer^n}c_k$
should converge weakly to a bona fide distribution in $\dist(\overline{\ddd}^n)$.  
In fact for any distributional solution $u$ to (\ref{sys}) the Fourier coefficients $c_k$ must be smooth functions, as follows. (In reference to lattice points $k\in\integer^n$ we use the convention established earlier whereby $k_{n+1}=0$.)   

\begin{prop}\label{prop-harmonic}
Any distributional solution $u$ to the equations (\ref{proof-wave}) has real analytic Fourier coefficients $c_k$ $(k\in\integer^n)$. Each such coefficient is a tensor product univariate functions, of the form
\begin{equation}\label{coefficient-structure}
c_k(w)=\prod_{j=1}^n\phi_j(w_j),\quad\mbox{ where }\quad-\ml_j\phi_j=k_jk_{j+1}\phi_j\quad(1\leq j\leq n),
\end{equation}
and is uniquely determined by its restriction to the set 
\[
\{w\in\overline{\ddd}^n\,|\,|w_n|=1\}.
\]
\end{prop}
\begin{pf}
If $u$ satisfies (\ref{proof-wave}) then for any $k\in\integer^n$ and any test function $\phi\in C^\infty(\overline{\ddd}^n)$
\[
\begin{split}
0&=\left[-\widetilde{\Delta}_ju+\frac{\partial^2u}{\partial\xi_{j+1}\partial\xi_j},\phi\otimes\mu_k\right]\\
&=[u,-\widetilde{\Delta}_j\phi\otimes\mu_k-k_jk_{j+1}\phi\otimes\mu_k]\\
&=[u,(-\widetilde{\Delta}_j\phi-k_jk_{j+1}\phi)\otimes\mu_k]\\
&=[c_k,-\widetilde{\Delta}_j\phi-k_jk_{j+1}\phi]\\
&=[-\widetilde{\Delta}_jc_k-k_jk_{j+1}c_k,\phi],
\end{split}
\]
which implies that 
\begin{equation}\label{coefficient-equation}
-\widetilde{\Delta}_jc_k-k_jk_{j+1}c_k=0.
\end{equation}
Since the operator $-\widetilde{\Delta}_j-k_jk_{j+1}$ is elliptic with real analytic coefficients it follows by analytic hypoellipticity \cite[Theorem~10]{Ha:2006} that each $c_k$ is real analytic.   
Thus the boundary condition (\ref{proof-boundary}) may be interpreted in the sense of ordinary functions, where 
\[
u(w,\one)=\sum_{k\in\integer^n}c_k(w).
\]
Equations (\ref{coefficient-equation}) for $1\leq j\leq n$ imply furthermore that each $c_k$ is a tensor product of eigenfunctions of the hybrid laplacian, 
\[
c_k(w)=\prod_{j=1}^n\phi_j(w_j),
\]
with $-\widetilde{\Delta}_j\phi_j=k_jk_{j+1}\phi_j$.   Note in particular that $k_{n+1}=0$.  Therefore $\phi_n$ is a harmonic function, determined by its restriction to the boundary circle $|w_n|=1$.   
\end{pf}

The real analytic functions $c_k$ may themselves be expanded as Fourier series with radial coefficients in the form
\begin{equation}\label{coefficient-expansion}
c_k(w)=\sum_{l\in\integer^n}d_{k,l}(r)e^{i\langle l,\theta\rangle},
\end{equation}
where $r=(|w_1|,\ldots,|w_n|)$.  Denote the left shift operator by a tilde, so that 
\[
\tilde{k}=(k_2,k_3,\ldots,k_n,0).
\]
\begin{prop}\label{prop-u-coupling}
Equations (\ref{proof-coupling}) imply that in the expansion (\ref{coefficient-expansion}) $d_{k,l}=0$ unless $l=k-\tilde{k}$.  
\end{prop}
\begin{pf}
For $w\in\overline{\ddd}^n$ write 
\[
r=(|w_1|,\ldots,|w_n|),\qquad v=\left(\frac{w_1}{|w_1|},\ldots,\frac{w_n}{|w_n|}\right),
\]
and let $P:\overline{\ddd}^n\rightarrow[0,1]^n\times\ttt^n$ denote the change of variables
\[
w\stackrel{P}{\mapsto}(r,v).
\]
The coefficient $d_{k,l}$ is defined in terms of $u$ by the formula
\begin{equation}\label{d-formula}
[d_{k,l},\rho]=\bigl[c_k,(\rho\otimes\mu_l)\circ P\bigr]=\bigl[u,((\rho\otimes\mu_l)\circ P)\otimes\mu_k\bigr]
\end{equation}
for test functions $\rho\in C^\infty([0,1]^n)$.  If $u$ satisfies the equation (\ref{proof-coupling}) it follows that 
\[
\begin{split}
0&=[C_ju,((\rho\otimes\mu_l)\circ P)\otimes\mu_k]\\
&=[u,-C_j((\rho\otimes\mu_l)\circ P)\otimes\mu_k]\\
&=[u,-(l_j-k_j+k_{j+1})((\rho\otimes\mu_l)\circ P)\otimes\mu_k]\\
&=[d_{k,l},-(l_j-k_j+k_{j+1})\rho]\\
&=[-(l_j-k_j+k_{j+1})d_{k,l},\rho],
\end{split}
\]
where 
\[
C_j=\frac{\partial }{\partial\theta_j}-\frac{\partial }{\partial\xi_j}+\frac{\partial }{\partial\xi_{j+1}}.
\]
Therefore $-(l_j-k_j+k_{j+1})d_{k,l}=0$ for $1\leq j\leq n$, so that $d_{k,l}=0$ unless $l=k-\tilde{k}$, as claimed.  \end{pf}

Relabeling appropriately, we may thus express the functions $c_k$ in the form 
\begin{equation}\label{d_k}
c_k(w)=d_k(r)e^{i\langle k-\tilde{k},\theta\rangle}=\prod_{j=1}^n\phi_j(w_j),
\end{equation}
consistent with (\ref{coefficient-structure}).  This yields a formal Fourier series for $u$ having smooth coefficients,
\begin{equation}\label{formal-fourier}
u(w,z)\sim\sum_{k\in\integer}d_k(r)e^{i\langle k-\tilde{k},\theta\rangle}z^k.
\end{equation}
The partial boundary condition (\ref{proof-boundary}) requires that for $w\in\partial\overline{\ddd}^n$, 
\begin{equation}\label{boundary-revised}
\sum_{k\in\integer}d_k(r)e^{i\langle k-\tilde{k},\theta\rangle}=\Psi(w,\one).
\end{equation}
In conjunction with Lemma~\ref{lem-complexification}, equation (\ref{d_k}) shows that if $k_jk_{j+1}\geq 1$, then $\phi_j$ has angular part $e^{i(k_j-k_{j+1})\theta_j}$ and is hence proportional to $\varphi^{(k_j,k_{j+1})}$.  Also, if $k_j>0$ and $k_{j+1}=0$, then $\phi_j$ is proportional to $\varphi^{(k_j,0)}$, since, up to a scalar multiple, the latter is the unique harmonic function with angular part $e^{ik_j\theta_j}$.  Lastly, if $k_jk_{j+1}<0$ for any $1\leq j\leq n$, then $\phi_j=0$ and $c_k=0$.   Thus $c_k\neq0$ only for lattice points $k$ for which $k_jk_{j+1}$ has constant sign.   

The structure of the function $\Psi(w,\one)$ further restricts the set of lattice points $k$ at which $c_k\neq0$ through condition (\ref{boundary-revised}).

Since $\Psi^{w_n}_1(0)=w_n$, the function $\Psi(w,\one)$ is holomorphic in  the variable $w_n$, and $\Psi(w,\one)$ is therefore determined by its restriction to the set 
\[
\{w\in\overline{\ddd}^n\,|\,|w_n|=1\}.
\]
Proposition~\ref{prop-harmonic} then implies that conditions (\ref{proof-boundary}) and (\ref{boundary-revised}) extend to all of $\overline{\ddd}^n$, leading to the following uniqueness result for the system (\ref{sys}).  
\begin{prop}\label{prop-uniqueness}
The system (\ref{sys}) has at most one distributional solution, necessarily a function 
\[
u\in L^2(\overline{\ddd}^n\times\ttt^n)
\]
of the form
\[
u(w,z)=\sum_{k\in\integer^n}d_k(r)e^{i\langle k-\tilde{k},\theta\rangle}z^k,
\]
where 
\[
d_k(r)=\frac{1}{(2\pi)^n}\int_{\ttt^n}\Psi\bigl(r_1e^{i\theta_1},\ldots,r_ne^{i\theta_n},\one\bigr)e^{-i\langle k-\tilde{k},\theta\rangle}\,d\theta.
\]
\end{prop}
\begin{pf}
Note that for every $w\in\overline{\ddd}^n$, $\Psi(w,\one)$ is a composition of disk automorphisms, evaluated at $0$, so that $|\Psi(w,\one)|\leq 1$.   
Since the coefficients $d_{k,l}$ are uniquely determined by the formula (\ref{d-formula}), the condition (\ref{boundary-revised})---extended to $\overline{\ddd}^n$---yields the given integral formula.   That $u\in L^2(\overline{\ddd}^n\times\ttt^n)$ then follows from the bound $|\Psi(w,\one)|\leq 1$.  In detail,
$(2\pi)^n\geq||\Psi(re^{i\cdot},\one)||^2_{L^2(\ttt^n)}=(2\pi)^n\sum_{k\in\integer}|d_k(r)|^2$.  Therefore,
\[
\begin{split}
||u(w,z)||^2_{L^2(\ddd^n\times\ttt^n)}&=\int_{\ddd^n\times\ttt^n}\bigl|\sum c_k(w)z^k\bigr|^2\,dwdz=(2\pi)^n\int_{\ddd^n}\sum|c_k(w)|^2\,dw\\
&=(2\pi)^{2n}\int_{[0,1]^n}\sum|d_k(r)|^2\,r_1\cdots r_n\,dr\leq 2^n\pi^{2n}.
\end{split}
\]
\end{pf}

The next step is to analyze the structure of $\Psi(w,\one)$ to determine the set of lattice points $k\in\integer^n$ for which $c_k\neq0$.   
\vspace*{5pt}
Let $\lll_n$ denote the set of all lattice points $k\in\integer^n$ with the properties: 
\begin{enumerate}
\item each $k_j\geq0$;
\item  $k_1=1$;
\item for each $2\leq j\leq n-1$, if $k_j=0$ them $k_{j+1}=0$.
\end{enumerate}  
For $w=\bigl(r_1e^{i\theta_1},\ldots,r_ne^{i\theta_n}\bigr)\in\overline{\ddd}^n$, let $g_k(r)$ denote the $k$th Fourier coefficient of 
\[
\Psi(w,\one)=\sum_{k\in\integer^n}g_k(r)e^{i\langle k,\theta\rangle}.
\]
\begin{prop}\label{prop-index-form}
If $g_l\neq0$ then there is a unique lattice point $k\in\lll_n$ such that $l=\widetilde{k}-k$.   
\end{prop}
\begin{pf}
If $l=k-\tilde{k}$ then the conditions that characterize $k\in\lll_n$ translate in terms of $l$ to:
\begin{enumerate}
\item for each $1\leq j\leq n$, $l_1+\cdots+l_j\leq 1$;
\item $l_1+\cdots+l_n=1$;
\item for each $1\leq j\leq n-1$, if $l_1+\cdots+l_j=1$ then $l_1+\cdots+l_{j+1}=1$.  
\end{enumerate}
As in \S\ref{sec-app-left} we argue by induction on the sequence of functions
\[
v=\Psi^{w_j}_{1}\circ\cdots\circ\Psi^{w_n}_{1}(0)
\]
as $j$ decreases from $n$ to $1$.  In the base case $v=\Psi^{w_n}_{1}(0)=w_n=e^{i\theta_n}$ there is a single one-dimensional vector, $1$, corresponding to a non-zero coefficient, and the associated set $\{1\}$ conforms the prescribed criteria.  For the induction, the key observation on the level of formulas is simply that 
\[
\begin{split}
\frac{w_j+v}{1+\overline{w}_jv}&=w_j+(1-r_j^2)\frac{v}{1+\overline{w}_jv}\\
&=w_j+(1-r_j^2)(v-\overline{w}_jv^2+\overline{w}_j^2v^3-\overline{w}_j^3v^4+\overline{w}_j^4v^5-\cdots).
\end{split}
\]
Each term $\overline{w}_j^sv^{s+1}$ corresponds to indices $(l_j,\ldots,l_n)$ where $l_j=-s$.  By inductive hypothesis we assume that conditions (1.-3.) are satisfied for the indices $(l^\prime_{j+1},\ldots,l^\prime_n)$ occurring in $v$.  It follows by the above formula that they are again satisfied for $(l_j,\ldots,l_n)$.  For example, the total sum being 1, means that the part of the sum coming from $v$ in a term $\overline{w}_j^sv^{s+1}$ is $s+1$, to which is added $-s$ from the term $l_j=-s$, for a total of $1$.  The other items are similar.   
\end{pf}

By equation (\ref{d_k}), Proposition~\ref{prop-index-form} implies that $c_k\neq0$ only if $k\in\lll_n$.   For such lattice points every function $\phi_j(w_j)$ in (\ref{d_k}) is proportional to a scattering polynomial, since only nonnegative indices $k_j$ occur, and the case $k_j=0, k_{j+1}>0$ is ruled out.  We thus have that the system (\ref{sys}) has at most one distributional solution, necessarily of the form
\begin{equation}\label{almost-structure}
u(w,z)=\sum_{k\in\lll_n}\beta_k\Bigl(\prod_{j=1}^n\varphi^{(k_j,k_{j+1})}(w_j)\Bigr)z^k,
\end{equation}
where each $\beta_k$ is a scalar given by the equation 
\begin{equation}\label{beta_k-computation}
\beta_k\prod_{j=1}^n\varphi^{(k_j,k_{j+1})}(w_j)=\frac{1}{(2\pi)^n}\int_{\ttt^n}\Psi\bigl(r_1e^{i\theta_1},\ldots,r_ne^{i\theta_n},\one\bigr)e^{-i\langle k-\tilde{k},\theta\rangle}\,d\theta.
\end{equation}
This completes the proof of Lemma~\ref{lem-unique}.

\end{document}